\documentclass[12pt]{amsart}
\usepackage{latexsym,amsfonts,amsthm,amsmath,amscd,amssymb}
\usepackage{hyperref}

\newtheorem{theorem}{Theorem}

\begin{document}

\title{Biharmonic space-like hypersurfaces in pseudo-Riemannian space}

\author{Wei Zhang}

\thanks{Mathematics Classification Primary(2000): 58E20.
\indent Thank Prof.Dong for the helpful instruction and Prof.Oniciuc
for kindly explanation of their work to me.}

\maketitle
\begin{abstract}
We classify the space-like biharmonic surfaces in 3-dimension
pseudo-Riemannian space form, and construct explicit examples of
proper biharmonic hypersurfaces in general ADS space.
\end{abstract}

\section{Introduction}
As a natural generalization of harmonic map, the biharmonic map was
introduced by Jiang who followed the idea of J.Eells and J.H.Sampson
in \cite{ES}, as the critical points of the bienergy
$$E_2(\phi)=\frac{1}{2}\int_M |\tau(\phi)|^2
v_g$$
where $\phi: (M,g)\rightarrow(N,h)$ is a smooth map between
two Riemannian manifolds and $\tau(\phi)$ is its first tension field
$trace \nabla d\phi$.

The biharmonic map also can be characterized by the vanishing of the
second tension field
$$\tau_2(\phi)=-\Delta^\phi \tau(\phi)-trace
R^N(d\phi,\tau(\phi))d\phi=0$$ and the non-harmonic one is called
proper biharmonic map.

Jiang had studied the biharmonic map in many general aspects(see
\cite{J1,J2}). While A.Balmus, S.Montaldo, C.Oniciuc, R.Caddeo and
E.Loubeau focused on the biharmonic isometric immersion, i.e. the
biharmonic submanifolds, and got some classification results(see
\cite{BB,BMO,CMO,MO} and their reference).

Pseudo-Riemannian space, especially the constant curvature ones: de
Sitter, Minkowski, anti de Sitter space with constant sectional
curvature 1, 0, -1, play important roles in the general relativity.

When the target manifold is pseudo-Riemannian manifold, the biengery
and second tension field can be defined in the same way, so are the
biharmonic map and biharmonic submanifold. Ouyang(\cite{Ou}) and
Sun(\cite{Su}) had study the space-like biharmonic submanifolds in
them. Remarkably, in \cite{Ou}, she had proved that
\begin{theorem}
$M^m$ is a biharmonic space-like hypersurface in de Sitter space
$S^{m+1}_1$ or Minkowski space $R^{n+1}_1$. If its mean curvatrue H
is constant, then it must be maximal.
\end{theorem}

Thus the  biharmonic submanifolds space-like in nonnegative
curvature space are more rigid. It is on the contrary to the
Riemannian case, for the generalized Chen conjecture says that {\it
biharmonic submanifolds of a manifold $N$ with $Riem^N \le 0$ are
minimal}.

In this paper, we first construct examples of proper biharmonic
hypersurfaces in ADS space.

Then we will classify the biharmonic space-like surfaces in
3-dimensional pseudo-Riemannian space form, namely

\begin{theorem}\label{Th:Main1}
If M is a proper biharmonic space-like surface in $H^3_1(1)$, then
it is $H^2(\frac{1}{\sqrt{2}})$ natrual embedded in.
\end{theorem}
and
\begin{theorem}\label{Th:Main2}
If M is a complete biharmonic space-like surface in $S^3_1$ or
$R^3_1$, then it must be totally geodesic, i.e. $S^2$ or $R^2$.
\end{theorem}

In effect, above two theorems can be viewed as the dual of their
Riemannian version.

\section{Preliminary}

\subsection{The biharmonic equation}

To study biharmonic submanifolds in space form, it is convenient to
split the bitension field in its normal and tangent components. In
fact, in Riemannian case, we have

\begin{theorem}[\cite {CMO}]
$M^m$ is a submanifold of $E^n(c)$, where $E^n(c)$ has constant
sectional curvature c, then it is biharmonic if and only if
\begin{equation}\label{Eq:BiharSub}
\begin{cases}
-\Delta^{\bot}H-traceB(\cdot, A_H \cdot)+mcH=0 &  \\
2traceA_{\nabla^{\bot}_{(\cdot)}H}(\cdot)+\frac{m}{2}grad(|H|^2)=0
\end{cases}
\end{equation}
where A denotes the Weingarten operator, B the second fundamental
form, H the mean curvature vector field, $\nabla^\bot$ and
$\Delta^\bot$ the connection and the Laplacian in the normal bundle
of M in $E^n(c)$
\end{theorem}

\subsection{The pseudo-Riemannian geometry}

Pseudo-Riemannian manifold is a manifold endowed with nondegenerated
but indefinite metric. Unlike the Riemannian metric, smooth manifold
may not carry pseudo-Riemannian metric with prescribed index, but if
one carries pseudo-Riemannian metric, there exists unique compatible
connection, the same with Riemannian geometry.

We focus on the pseudo-Riemannian space form with index one. From
now on, always denote $E^{m+1}_1(c)$ the $(m+1)$-dimension constant
curved pseudo-Riemannian space with index 1.

Except $|H|^2$ is a minus one, the biharmonic equation in
pseudo-Riemannian case coincides with Eq.~(\ref{Eq:BiharSub}) for
the connection still compatible with the indefinite metric.

If $M^m$ is a space-like hypersurface in $E^{m+1}_1(c)$, its
biharmonic equation can be simplified to(\cite{Ou})

\begin{equation}\label{Eq:BiharHyp}
\begin{cases}
-\Delta^{\bot}H-|B|^2H+mcH=0 &  \\
2traceA_{\nabla^{\bot}_{(\cdot)}H}(\cdot)+\frac{m}{2}grad(|H|^2)=0
\end{cases}
\end{equation}
It should be pointed out that the length square of the second
fundamental form $|B|^2$ and $|H|^2$ are both negative because $H$
has minus length.

\section{Example of proper biharmonic space-like hypersurfaces in the ADS space}

Some how, this section can be viewed as the dual case of
constructing biharmonic hypersurface in $S^n$, reader could refer to
the author's another paper \cite{Zh}.

The pseudohyperbolic space of radius r($>0$) is the hyperquadric
$$H^{n+p}_p(r)=\{x\in
R^{n+p+1}_{p+1};<x,x>=x_1^2+\cdots+x_n^2-x_{n+1}^2-\cdots-x_{n+p+1}^2=-r^2\}$$
this is a space of constant curvature $-\frac{1}{r^2}$, with index
$p$. Let $H^n(r)$ be the component of $H^n_0(r)$ through $(0,
\ldots, 0, r)$. It is the usual hyperbolic space.

\subsection{First example: Hyperbolic space}
Embedding $H^{n+1}_1(1)$ in $R^{n+2}_2$ and  fixing $x_{n+2}=s$,
$|s|<1$, we get a space-like submanifold $H^n(\sqrt{1-s^2})$. The
same with the hypersphere in the $S^{n+1}$, it has constant mean
curvature and $|B|^2=n\frac{s^2}{s^2-1}$. By the biharmonic equation
(\ref{Eq:BiharHyp}), we know $H^n(\frac{1}{\sqrt{2}})$ is a proper
biharmonic space-like hypersurface in ADS space.

\subsection{Second example: Clifford torus}

 Let $u$, $v$ be points of $H^m(r)$ and $H^n(s)$
respectively, where $r^2+s^2=1$. Then $x=(u,v)$ is a vector in
$R^{m+n+2}_2$ with $<x,x>=-1$.

This defines an isometric immersion of $H^m(r) \times H^n(s)$ into
$H^{m+n+1}_1(1)$, i.e, $H^m(r) \times H^n(s)$ is a space-like
hypersurface of the anti de Sitter space.

Obviously, it has constant mean curvature. The next thing we have to
do is to manipulate the radius to make sure it satisfying the
biharmonic equation (\ref{Eq:BiharHyp}).

Use the moving frame method, the second fundamental form of $H^m(r)
\times H^n(s)$ is
$$
\left(
  \begin{array}{ccccccc}
    \frac{s}{r} &  &  &  &  &  & 0 \\
     & \ddots &  &  &  &  &  \\
     &  & \frac{s}{r} &  &  &  &  \\
     &  &  & \ddots &  &  &  \\
     &  &  &  & \frac{r}{s} &  &  \\
     &  &  &  &  & \ddots &  \\
    0 &  &  &  &  &  & \frac{r}{s} \\
  \end{array}
\right)
$$
its length square is $-(\frac{s^2}{r^2}m+\frac{r^2}{s^2}n)$. Denote
$t=\frac{s^2}{r^2}$, the biharmonic equation becomes a quadric
equation about $t$. Solve it and exclude the maximal one, we have
$t=1$.

When $m \neq n$, $H^m(\frac{1}{\sqrt{2}}) \times
H^n(\frac{1}{\sqrt{2}})$ is a proper biharmonic space-like
hypersurface in ADS space.

\section{classification of biharmonic space-like hypersurface in 3-dimension $E^3_1(c)$}
Follow the idea of Chen and Ishikawa (\cite{CI}) or Caddeo
et.al.(\cite{CMO}, denote $\eta$ the unit vector of the normal
bundle. It is a time-like vector field, i.e $|\eta|=-1$. Then the
mean curvature can be written as $H=f\eta$, the biharmonic equation
(\ref{Eq:BiharHyp}) becomes
\begin{equation}
\begin{cases}\label{Pse:Eq:BiharSur}
\Delta f=(|A|^2+2c)f &\\
A(grad f)-f grad f=0 &
\end{cases}
\end{equation}
where the Laplacian is in Beltremi sense, minus of the trace
Laplacian and $|A|^2$ is a positive norm.

By the same way of proving theorem 4.5 in \cite{CMO}, we can show
that $f$ must be constant. For readers' convenience, we still write
down the details.

Using the contradiction method, suppose that $grad f$ do not
vanishing in a neighborhood $U$, let $X_1$, $X_2$ be a orthpnormal
frame, where $X_1=\frac{grad f}{|grad f|}$. Then $X_2 f=0$, $X_1
f=|grad f|$.

Then the second fundamental form $B$ of $M$ is given by
\begin{equation}\label{Eq:SecondFundf}
 B(X_1,X_1)=-f\eta, \quad B(X_1,X_2)=0, \quad B(X_2,X_2)=3f\eta,
\end{equation}

Therefore, $|A|^2=10f^2$.

Take the derivative of the second fundamental form
(\ref{Eq:SecondFundf}), by the Codazzi equation,
\begin{equation}\label{Pse:Eq:Connec}
X_2f=-4f\omega_{21}(X_1); 3X_1f=-4f\omega_{12}(X_2)
\end{equation}
Lead us to $\omega_{21}(X_1)=0$, moreover, $d\omega_1=0$.

Thus $\omega_1$ is locally exact. We can find function $u$, s.t.
$du=\omega_1$. At the same time, for $df=(X_1 f)\omega_1+(X_2
f)\omega_2$ and $X_2 f=0$, $df\wedge \omega_1=0$, i.e. $f$ is a
function only about $u$. Denoting $f$ as $f(u)$, the first and
second derivative of $f$ about $u$ are denoted as $f'$ and $f''$.

By the second equation of (\ref{Pse:Eq:Connec}), there is:
$$4f\omega_{12}=-3f'\omega_2$$

Based on (\ref{Pse:Eq:Connec}), the second covariant derivative and
Laplace of $f$ can be calculated,
$$4f\Delta f=4f(-tr(Hess(f)))=3(f')^2-4ff''$$

By the first equation of (\ref{Pse:Eq:BiharSur}),
\begin{equation}\label{Pse:Eq:Codaf}
3(f')^2-4ff''=40f^4+8cf^2
\end{equation}
Let $y=(f')^2$, it becomes:
$$2f\frac{dy}{df}-3y=-40f^4-8cf^2$$

Solve it, we have:

\begin{equation}\label{Pse:Eq:CodaPn}
(f')^2=-8f^4-8cf^2+Cf^{\frac{3}{2}}
\end{equation}

At another hand, using Gauss equation, we get
$$K=c-detA$$
The minus sign before $detA$ is because of the normal vector is
time-like and having negative length.

Combining the structure equation, there are
\begin{equation}
\begin{cases}
K=c+3f^2\\
d\omega_{12}=-K\omega_1 \wedge \omega_2
\end{cases}
\end{equation}

Followed by
\begin{equation}\label{Pse:Eq:Gausf}
4ff''-7(f')^2-16f^4-\frac{16c}{3}f^2=0
\end{equation}

Comparing (\ref{Pse:Eq:Codaf}) and (\ref{Pse:Eq:Gausf}), there is
\begin{equation}\label{Pse:Eq:GausPn}
(f')^2=-14f^4-\frac{10c}{3}f^2
\end{equation}

Combine (\ref{Pse:Eq:CodaPn}) and (\ref{Pse:Eq:GausPn}), no matter
 what $c$ is, $f$ must be constant, contradict to our assumption.

In the case of $c=-1$, we have $|A|^2=2$. Denoting $-\lambda_1,
-\lambda_2$ the eigenvalue of the Weingarten operator, then them
must be constant. Using the Gauss and Codazzi equation, we could
conclude that $\lambda_1=\lambda_2=\pm 1$ or
$\lambda_1=-\lambda_2=\pm 1$.

To see this, still denote $\{X_1, X_2\}$ the frame diagonalizing
$A$, and $\omega_{ij}$ the connection form respectively, then
\begin{equation}\label{Eq:SecondFund}
 B(X_1,X_1)=\lambda_1\eta, \quad B(X_1,X_2)=0, \quad B(X_2,X_2)=\lambda_2\eta,
\end{equation}

Take the derivative $B(X_2,X_2)$ in $X_1$ and notice
equation(\ref{Eq:SecondFund}),

\begin{equation*}
\begin{split}
&0=\nabla_{X_1}(B(X_2,X_2))\\
=&\nabla_{X_2}(B(X_1,X_2))-B(\nabla_{X_2}X_1,X_2)-B(X_1,\nabla_{X_2}X_2)+2B(\nabla_{X_1}X_2,X_2)\\
=&-\omega_{12}(X_2)\lambda_2\eta-\omega_{21}(X_2)\lambda_1\eta
\end{split}
\end{equation*}

case 1, $\omega_{12} \neq 0$, then $\lambda_1=\lambda_2=\pm 1$;

case 2, $\omega_{12}=0$, the surface is flat. By the Gauss equation,
there is

$$-1=-(-\lambda_1)\lambda_2$$
then $\lambda_1=-\lambda_2=\pm 1$.

As the fundamental theorem of surface can be easily generalized to
constant curved space forms, the first one corresponding to
hyperbolic space in our example and the second the Clifford torus.
For the biharmonic Clifford torus of dimension 2 must be maximal,
the only proper biharmonic surface is the $H^2(\frac{1}{\sqrt{2}})$.
Finish the proof of theorem \ref{Th:Main1}.

We could get more rigidity in the Nonnegative curved case. By the
same method, we could show any biharmonic space-like hypersurface in
3-dimensional de Sitter or Minkowski space has constant mean
curvature. According to Ouyang's results, it must be maximal.
Finally, combining Ishihara (\cite{Is}, theorem 1.1) and Calabi's
(\cite{Ca}, or refer to \cite{CY}, page 417, corolary) results
respectively, it must be geodesic. Theorem \ref{Th:Main2} follows.

This section is somehow an analog of classification of biharmonic
surface in 3-dimensional Riemannian space form(\cite{CMO}).

\vfill

\noindent Wei Zhang

\

\noindent School of Mathematical Sciences

\noindent Fudan University

\noindent Shanghai, 200433, P. R.China

\

\noindent Email address: 032018009@fudan.edu.cn

\end{document}